\documentclass{amsart}
\usepackage{amssymb, amsmath}
\usepackage{fullpage}
\usepackage{array}
\usepackage{mathrsfs}
\usepackage{xcolor}
\usepackage{hyperref}
\usepackage{graphicx,color}
\usepackage{float}
\usepackage{tikz}
\usepackage{comment}
\usetikzlibrary{decorations.pathreplacing,calligraphy}
\usepackage{}

\newcommand{\co}[1]{\textcolor{blue}{#1}} 

\newcommand{\R}{\mathbb{R}}
\newcommand{\N}{\mathbb{N}}

\newcommand{\Z}{\mathbb{Z}}

\DeclareMathOperator{\Rn}{R_{\beta \mbox{,} k}}

\DeclareMathOperator{\ase}{A_{\sigma \eta}}
\DeclareMathOperator{\SL(2)}{SL(2,\R)}

\DeclareMathOperator{\ta}{\tau_{Z}}
\DeclareMathOperator{\cx}{C_1^{q_1}}
\DeclareMathOperator{\cy}{C_2^{q_2}}	
\newtheorem{teorema}{Theorem}[section] 
\newtheorem{exemplo}[teorema]{Example}

\newtheorem{lema}[teorema]{Lemma}
\newtheorem{corolario}[teorema]{Corollary}
\newtheorem{obs}[teorema]{Remark}

\newtheorem*{TEOA}{Theorem A}

\newcolumntype{C}[1]{>{\centering\let\newline\\\arraysetminus\hspace{0pt}}m{#1}}

\usepackage[utf8]{inputenc}
\usepackage{amsfonts}
\usepackage{amssymb}
\usepackage{verbatim}
\numberwithin{equation}{section}

\begin{document}
\title{Discontinuity of Lyapunov exponents for $\SL(2)$-valued cocycles}

\author{Edhin Mamani$^1$}
\address{ICEx, UFMG, Avenida Antonio Carlos 6627,
31123-970, Belo Horizonte, Brazil}
\email{emamani@ufmg.br}

\author{Raquel Saraiva$^{2}$}
\address{ICEx, UFMG, Avenida Antonio Carlos 6627,
31123-970, Belo Horizonte, Brazil}
\email{rss07@ufmg.br}

\thanks{\footnotesize $^{1}$ Partially supported by FAPEMIG and CNPq}
\thanks{\footnotesize $^{2}$ Supported by FAPEMIG}


\begin{abstract}
We exhibit an example of discontinuity point for the Lyapunov exponents as a function of the cocycle in the $\alpha$-Hölder topology. The linear cocycle taking values in $\SL(2)$ is locally constant and defined over a Bernoulli shift. Our example extends Bocker-Viana's and Butler's results. In particular, it gives a partial answer to a question raised by Clark Butler. Finally, we give an example of discontinuity in the setting of $\mathrm{SL}(m,\R)$-valued cocycles, which is constructed from $\SL(2)$-valued cocycles.
\end{abstract}

\maketitle

\section{Introduction}
Lyapunov exponents, first introduced by Lyapunov in 1892 to study the stability of solutions to ordinary differential equations \cite{lya92}, have since become a fundamental tool in the analysis of smooth dynamical systems, particularly for $C^k$-diffeomorphisms on manifolds. For a $C^{1+\alpha}$-diffeomorphism preserving an ergodic measure with one positive and one negative Lyapunov exponent, Pesin (1976) established the almost-everywhere existence of local stable and unstable manifolds \cite{pes76}. He later generalized this theory to measures with all non-zero Lyapunov exponents, introducing what is now known as Pesin theory, which extends key aspects of uniformly hyperbolic (Axiom A) systems to non-uniformly hyperbolic settings \cite{pes77}. Beyond hyperbolicity, Lyapunov exponents are deeply connected to metric entropy and the fractal dimension of invariant measures, as demonstrated in works by Barreira, Ledrappier, and Young \cite{bar12,led86,you82}. Due to the challenges inherent in smooth dynamics, the concept has been extended to non-smooth settings, where they are studied via linear cocycles over continuous base dynamical systems.

An important question in the theory of Lyapunov exponents concerns their regularity properties. From Oseledets' work, we know that under certain hypotheses, Lyapunov exponents vary measurably with respect to both the cocycle and the invariant measure \cite{oseledets}. However, the most significant regularity property is the continuity of Lyapunov exponents with respect to variations in the linear cocycles. This continuity property is particularly valuable as it serves as a powerful tool for deriving results about statistical and topological properties of dynamical systems. For instance, continuous variation of exponents enables the study of stability properties and bifurcation phenomena \cite{viana}, and plays a crucial role in the analysis of random matrix products \cite{bou16}. Furthermore, for certain systems continuity is an important tool to prove the strong positive recurrence property defined by Buzzi-Crovisier-Sarig, which in turn implies several ergodic consequences for the system \cite{bcs22,bcs25}. 

We study linear cocycles over continuous dynamical systems, with particular focus on $\mathrm{SL}(2,\mathbb{R})$-valued cocycles. Consider a continuous map $A \colon M \to \mathrm{SL}(2, \mathbb{R})$ over a compact metric space $(M,d)$. The linear cocycle defined by $A$ over the base dynamics $f:M \to M$ is the skew-product $F_A \colon M \times \mathbb{R}^2 \to M \times \mathbb{R}^2$ defined by $F_A(x, v) = (f(x), A(x)v)$. Since the base transformation $f$ is fixed, we always identify the cocycle $F_A$ with its generator $A$. 

The cocycle property $A^{n+m}(x)=A^n(f^m(x))A^m(x)$ for all $n,m\in \Z$ follows directly from the dynamical definition of the iterated products: 
\begin{equation}\label{lll}
    A^n(x) = 
\begin{cases}
A(f^{n-1}(x)) \cdots A(x) & \text{if } n > 0, \\
\text{Id} & \text{if } n = 0, \\
A(f^{-n}(x))^{-1} \cdots A(f^{-1}(x))^{-1} & \text{if } n < 0.
\end{cases}
\end{equation}

For every $f$-invariant probability measure $\mu$ on $M$, the Lyapunov exponents associated to the system $(f,A,\mu)$ are defined by
\[   \lambda_{+}(A,\mu,x) = \lim_{n\rightarrow \infty}\frac{1}{n}\log\|A^n(x)\| \quad \text{ and } \quad \lambda_{-}(A,\mu,x) = \lim_{n\rightarrow \infty}\frac{1}{n}\log\|A^n(x)^{-1}\|^{-1} \]
where the Furstenberg-Kesten Theorem guarantees the existence of these limits $\mu$-almost everywhere \cite{furstenberg}. Since Lyapunov exponents are $f$-invariant functions, they are constant $\mu$-almost everywhere when $\mu$ is ergodic. Therefore, for an ergodic measure $\mu$, we simply denote the exponents by $\lambda_+(A,\mu)$ and $\lambda_-(A,\mu)$. Moreover, $\lambda_+(A,\mu)$ and $\lambda_-(A,\mu)$ are upper semi-continuous and lower semi-continuous functions of $A$ respectively. For $\mathrm{SL}(2,\R)$-valued cocycles, if $\lambda_{+}(A, \mu)$ is zero then so is $\lambda_{-}(A, \mu)$ since $\lambda_{+}(A, \mu)+\lambda_{-}(A, \mu)=0$, and hence $A$ is a continuity point.

Fixing an ergodic measure $\mu$, Mañé-Bochi's Theorem provides the first information about continuity of Lyapunov exponents of an aperiodic system $(f,A,\mu)$ \cite{bochi}. This theorem establishes that, in the space of continuous cocycles $C^0(M,\mathrm{SL}(2,\mathbb{R}))$ over $(M,f,\mu)$, the only continuity points for Lyapunov exponents are either uniformly hyperbolic cocycles or cocycles with zero Lyapunov exponents. In particular, discontinuity is a generic phenomenon: any typical cocycle with nonzero exponents can be approximated in the $C^0$ topology by cocycles with zero Lyapunov exponents. Previously in the 1980s, Mañé also stated a similar conclusion for derivative cocycles of area-preserving $C^1$ diffeomorphisms on surfaces \cite{mane}. 

While discontinuities of Lyapunov exponents are common, there are notable situations for which continuity has been verified. For random products of $2 \times 2$ matrices, continuity was proved in the Bernoulli and Markov settings by Bocker-Viana~\cite{bockerviana} and Malheiro-Viana~\cite{malheiroviana} respectively. In the higher dimensional context, the picture is more subtle, but recently Avila, Eskin and Viana~\cite{aev} announced continuity for i.i.d. random products of $m\times m$ matrices, indicating that such stability is not restricted to the low-dimensional setting.

For $l\geq1$, let $M = \{0,\cdots, l\}^{\mathbb{Z}}$ be the space of bi-infinite sequences of $l+1$ symbols and $f \colon M \to M$ be the left shift map defined by $f\big((x_i)_{i \in \mathbb{Z}}\big) = (x_{i+1})_{i \in \mathbb{Z}}$. We equip $M$ with the metric 
\begin{equation}\label{metric}  
d(x,y) = 2^{-N(x,y)} \quad \text{ where } \quad N(x,y) = \min\{|i| \geq 0 : x_i \neq y_i\}, 
\end{equation}
which in turn generates the product topology on $M$. Given $p_i \in (0,1)$ for $i\in \{0,\cdots,l\}$, consider the probability measure $\nu = \sum_{i=0}^lp_i\delta_i$ on $\{0,1\}$ and let $\mu = \nu^{\mathbb{Z}}$ be the corresponding Bernoulli product measure on $M$. This construction yields the Bernoulli shift system $(M, f, \mu)$.

One may consider coarser topologies than the well-studied $C^0$ topology. Given $\alpha > 0$, the $\alpha$-Hölder topology is induced by the norm
\begin{equation}\label{norma}
\|A\|_{\alpha} = \|A\|_{\infty} + \sup_{x \neq y} \frac{\|A(x)-A(y)\|}{d(x,y)^{\alpha}},
\end{equation}
where $\|A\|_{\infty} = \sup_{x\in M}\|A(x)\|$. Denote by $C^{\alpha}(M,\mathrm{SL}(2,\mathbb{R}))$ the space of $\alpha$-Hölder continuous cocycles equipped with this topology. In the $\alpha$-Hölder topology, Backes-Brown-Butler \cite{bbb} established continuous variation of Lyapunov exponents for cocycles $A \in C^\alpha(M,\mathrm{SL}(2,\mathbb{R}))$ satisfying the $\alpha$-fiber-bunching condition: there exists an integer $N \geq 1$ such that for all $x \in M$,
\[  \|A^N(x)\| \cdot \|A^N(x)^{-1}\| < 2^{\alpha N}.  \]
Within the class of $\alpha$-Hölder continuous cocycles that fail to satisfy the fiber-bunching condition, Bocker and Viana exhibited explicit examples of discontinuity points for the Lyapunov exponents \cite{bockerviana}. More precisely, for each parameter $\sigma > 1$, they considered the cocycle $A_\sigma \colon \{0,1\}^{\Z} \to \mathrm{SL}(2,\mathbb{R})$ defined by
\begin{equation}\label{cla}
    A_\sigma(x) = 
\begin{cases}
    \begin{pmatrix}
        \sigma^{-1} & 0 \\
        0 & \sigma 
    \end{pmatrix} & \text{if } x_0 = 0, \\[2ex]
    \begin{pmatrix}
    \sigma & 0 \\
    0 & \sigma^{-1}
    \end{pmatrix} & \text{if } x_0 = 1.
\end{cases}
\end{equation}
The Lyapunov exponents of the cocycle $A_{\sigma}$ can be computed explicitly via the Birkhoff Ergodic Theorem, yielding 
\begin{equation}\label{lyapsigma}
\lambda_{\pm}(A_{\sigma}, \mu)=\pm|(1-2p)\log \sigma| \quad \text{ where } \quad p=\nu(\{1\}).
\end{equation}
This cocycle exhibits nonzero Lyapunov exponents for all Bernoulli measures $\nu$ induced by $p \neq 1/2$. Furthermore, $A_{\sigma}$ satisfies the $\alpha$-fiber-bunching condition precisely when $\sigma^2 < 2^\alpha$, making it a continuity point for the Lyapunov exponents according to the Backes-Brown-Butler's Theorem \cite{bbb}. On the other hand, Bocker-Viana's construction specifically applies to cocycle $A_\sigma$ with parameters satisfying
\begin{equation}\label{viva}
    \sigma^2 > 2^{4\alpha} \quad \text{and} \quad p\neq \tfrac{1}{2}.
\end{equation}
By refining the techniques of Mañé-Bochi and Knill \cite{knill}, Bocker and Viana showed that $A_\sigma$ can be approximated in the $\alpha$-Hölder topology with cocycles having zero Lyapunov exponents. Notably, observe that these Bocker-Viana cocycles are strongly non-fiber-bunched. Building on this technique, Butler \cite{butler} constructed cocycles that are arbitrarily close to satisfying the $\alpha$-fiber-bunching condition while still serving as discontinuity points for the Lyapunov exponents in the $\alpha$-Hölder topology. Specifically, Butler proved this for cocycles $A_\sigma$ with parameters satisfying
\begin{equation}\label{mejor}  \sigma^{4p-2} \geq 2^{\alpha} \quad \text{ and } \quad p \in \left(\tfrac{1}{2},1\right).   \end{equation}
These discontinuity points approach the $\alpha$-fiber-bunched regime as $p \to 1^-$ for the Bernoulli measure $\nu$ induced by $p$. A crucial feature of Butler's construction is that approximating cocycles now exhibit arbitrarily small but nonzero Lyapunov exponents, which stands in direct contrast to the zero Lyapunov exponents in earlier results. 

In \cite{butler}, Butler also posed the question of whether approximation by cocycles with zero Lyapunov exponents remains possible when $\sigma^2 \in (2^{2\alpha}, 2^{4\alpha})$. As a consequence of our main theorem, we provide a partial answer to this question.
\begin{teorema}\label{coro}
For any $\alpha>0$ and $\sigma > 1$ satisfying $\sigma^2 \geq 2^{3\alpha}$ there exist $\alpha$-Hölder continuous cocycles $B \colon \{0,1\}^{\Z} \to \mathrm{SL}(2,\mathbb{R})$ with zero Lyapunov exponents that are arbitrarily close to $A_\sigma$ in the $\alpha$-Hölder topology. In particular, $A_\sigma$ is a discontinuity point for the Lyapunov exponents with respect to $\mu$ in the space $C^\alpha(\{0,1\}^{\Z},\mathrm{SL}(2,\mathbb{R}))$.
\end{teorema}
This result follows from an argument based on interchanging the Oseledets subspaces of $A_{\sigma}$ (see Section 3 for details). In contrast, the proof of zero Lyapunov exponents uses a new argument that relies on the explicit structure of the perturbation. This provides an original approach to showing the cancellation of Lyapunov exponents. While this theorem proves discontinuity in the stated parameter range, our approach does not extend to the more general setting described in \eqref{mejor}. The known results for the cocycle $A_\sigma$ can be summarized in the following diagram.

\bigskip

\begin{tikzpicture}[decoration=brace]
\draw(0,0)--(13,0);
\foreach \x/\xtext in {0/$0$,1/$\sigma^2$,3/$2^{\alpha}$,5/$2^{2\alpha}$,7/$\sigma^2$,9/$2^{3\alpha}$,10/$\sigma^2$,13/}
  \draw(\x,1pt)--(\x,-1pt) node[above] {\xtext};
\draw[decorate, line width=0.3mm, decoration={calligraphic brace, amplitude=0.1cm, mirror}, yshift=-2ex]  (0,0) -- node[below=0.4ex] {$\begin{array}{lll}
\hspace{1.3cm}\mbox{Continuity} \\
	\hspace{1.4cm}\forall p \in [0,1]\\
 \hspace{0.8cm}\mbox{(Backes, Brown, Butler)}
 \end{array}$}  (3,0);
\draw[decorate, line width=0.3mm, decoration={calligraphic brace, amplitude=0.1cm, mirror},yshift=-2ex]  (5,0) -- node[below=0.4ex] {$\begin{array}{llll}
\hspace{1.5cm}\mbox{Discontinuity for} \\
\hspace{1.8cm}p \in (3/4,1)  \\
\hspace{0.6cm} \mbox{(Consequences of Butler's result)}
 \end{array}$}  (8.89999,0);
\draw[decorate,line width=0.3mm, decoration={calligraphic brace, amplitude=0.1cm, mirror},  yshift=-2ex]  (9.011111,0) -- node[below=0.4ex] {$\begin{array}{lll}
\hspace{0.8cm}\mbox{Discontinuity for} \\
\hspace{1.5cm}p \neq \frac{1}{2} \\ \hspace{0.8cm} \mbox{(Mamani-Saraiva)}
\hspace{0cm} 
 \end{array}$}(13,0);
 \end{tikzpicture}
 
\bigskip

Given parameters $a_i>0$ with $i\in \{0,\cdots,l\}$ such that there exist two indices $r$ and $j$ satisfying $a_r>1$ and $a_j<1$, we define a family of cocycles $A: M \to \mathrm{SL}(2,\mathbb{R})$ by 
\begin{eqnarray}\label{new}
A|_{[0;i]}=\begin{pmatrix}
				a_{i} & 0 \\
				0 & a_{i}^{-1} 
		\end{pmatrix}\end{eqnarray}

This family generalizes cocycle $A_\sigma$ defined in \eqref{cla}, making $l=1$ and $a_r= a_j^{-1}=\sigma$. The Lyapunov exponents of cocycle $A$ can be computed explicitly as follows
\begin{equation}\label{lyapgeneral}\lambda_{\pm}(A, \mu)=\pm\left|\sum_{a_k>1}p_k\log a_k-\sum_{a_k<1}p_k\log a_k^{-1}\right|
\end{equation}

We also require that weights $(p_0, \ldots, p_l)$ be chosen so that $\lambda_+(A, \mu) \neq 0$.
Moreover, since $A$ depends only on the first coordinate it follows that $A$ is $\alpha$-Hölder continuous for every $\alpha>0$.

For the sake of simplicity, we define 
\begin{eqnarray}\label{min}
    \eta:= \min \{a_j^{-1}, a_r\} \,\ \,\ \mbox{ and } \,\ \,\ \sigma:=\max \{a_j^{-1}, a_r\}.
\end{eqnarray}

The following result represents the main contribution of our work and provides another concrete and explicit example of discontinuity in the $\alpha$-Hölder topology for the Lyapunov exponents. 

\begin{TEOA}\label{raquel}
\textit{For any $\alpha>0$, with $a_j<1<a_r$ and positive weights $(p_0,\ldots, p_l)$ satisfying $\lambda_+(A, \mu)\neq 0$ such that  
    \begin{eqnarray}\label{cond2}
            2^{3\alpha}<\sigma^2
        \end{eqnarray}
        and
\begin{eqnarray}\label{cond1}
2^{\left(2+\frac{\log\eta}{\log \sigma}\right)\alpha}\leq \eta^2
\end{eqnarray}
  there exist $\alpha$-Hölder continuous cocycles $B: M \to \SL(2)$ with $\lambda_{\pm}(B, \mu)=0$ which are arbitrarily close to $A$ in the $\alpha$-Hölder norm. Consequently, $A$ is a discontinuity point for Lyapunov exponents in $C^{\alpha}(M, \SL(2))$.}
	\end{TEOA}

\begin{obs}
    Condition \eqref{cond2} ensures that inequality \eqref{cond1} is not empty, that is, under this assumption there exists a constant $\eta<\sigma$ satisfying \eqref{cond1}.
\end{obs}
Note that by setting $\sigma = \eta$, we recover Theorem~\ref{coro}. 
Furthermore, hypothesis of Theorem A involve only $\eta$ and $\sigma$, which are defined in \eqref{min}. Therefore, no additional assumptions are needed for parameters $a_i$ in definition of $A$ when $a_i<\eta$. For instance, if either (i) $a_k<1$ and $a_k^{-1}<\eta$, or (ii) $a_k>1$ and $a_k<\eta$, then there are no restrictions required by parameters $a_k$. Thus, the construction remains valid whether $a_k$ satisfies the bound $2^{\left(2+\frac{\log\eta}{\log \sigma}\right)\alpha}$ or not. This reveals a key feature of our framework: the argument is decisively governed by the largest parameter $\eta$ leaving the others largely unconstrained. To illustrate a possible configuration, consider choosing $a_k<1$ such that $2^{\alpha}<a_k^{-2}<2^{\left(2+\frac{\log \eta}{\log \sigma}\right)\alpha}$.

\bigskip

 \begin{tikzpicture}[decoration=brace]
\draw(0,0)--(13,0);
\foreach \x/\xtext in {0/$0$,1/$a_s$, 1.8/$\eta^{-1}$,2.4/$a_i$,2.9/$a_k$,3.5/1,5/$2^{\alpha}$,6/$a_k^{-2}$,7.4/$2^{\left(2+\frac{\log \eta}{\log \sigma}\right)\alpha}$,9/$a_i^{-2}$,10.2/$\eta^2$,11/$\sigma^2$,12.4/$a_s^{-2}$}
  \draw(\x,1pt)--(\x,-1pt) node[above] {\xtext};
\draw[decorate, line width=0.3mm, decoration={calligraphic brace, amplitude=0.1cm, mirror},yshift=-2ex]  (7.4,0) -- node[below=0.4ex] {$\begin{array}{llll}
\hspace{0.7cm}\mbox{Discontinuity for} \\
\hspace{1.1cm} \mbox{Theorem~A} 
 \end{array}$}  (13,0);
\end{tikzpicture}

The condition in Theorem A depends on the smallest constant because of the underlying construction of the perturbation, which is supported in a certain cylinder. Specifically, the cylinder is built by first placing the symbols corresponding to $\eta$ (which can be $j$ or $r$) and then those for $\sigma$, with variation arising solely in the frequency of each symbol. Since $\eta$ is the smaller constant, its matrix produces weaker expansion in the vertical direction (if $\eta=a_r$) or weaker contraction in the horizontal direction (if $\eta=a_j^{-1}$) compared to the matrix for $\sigma$. This disparity means the first exchange requires more iterates than the second. The resulting reduction in iterates for the second exchange forces its closing angle to be smaller, which is the direct source of our condition. 

Finally, we provide a simple example of discontinuity of Lyapunov exponents in the $\alpha$-Hölder topology for $\mathrm{SL}(m,\R)$-valued cocycles. The cocycles considered are constructed with block diagonal matrices where each block is the $2\times 2$ $A_\sigma$ (or $A_\sigma^{-1}$) of Equation \eqref{cla} (See Section \ref{calda} for details). Thus, we obtain an analogous example in higher dimension to that given in Theorem \ref{coro}.
\begin{corolario}\label{dim}
For any $\alpha>0$, $\sigma > 1$ and $m\in\N$ satisfying $\sigma^2 > 2^{3\alpha}$ there exist cocycles $A \colon \{0,1\}^{\Z} \to \mathrm{SL}(2m,\mathbb{R})$ and $\tilde{A} \colon \{0,1\}^{\Z} \to \mathrm{SL}(2m+1,\mathbb{R})$, which are discontinuity points for $2m$ Lyapunov exponents in the $\alpha$-Hölder topology.    
\end{corolario}
Observe that for the even case $2m$ we obtain cocycles $A$ that are discontinuity points for all possible Lyapunov exponents. For the odd case, however, the cocycles $\tilde{A}$ are discontinuity points only for all Lyapunov exponents except a single one.

\section{Construction}\label{section2}
This section is devoted to half of the proof of Theorem~A, specifically to the construction of a suitable cocycle~$B_k$ that is arbitrarily close to~$A$ in the $\alpha$-H\"older topology. This construction refines the technique developed by Bocker and Viana in~\cite{bockerviana}, which, in turn, builds on methods from Mañé-Bochi, \cite{bochi}, and Knill, \cite{knill}. We postpone the other half of the proof of Theorem~A, the proof of zero exponents, to the next section. As mentioned in the introduction, we establish the vanishing of the Lyapunov exponent by directly exploiting the definition of the perturbation and the fact that the cocycle is constant on cylinders.

For the remainder of the proof, we assume $\sigma = a_r$ and $\eta = a_j^{-1}$. The other case is analogous but requires a slightly different perturbation, which is composed of a rotation matrix and two shear matrices. These shear matrices preserve the horizontal direction but not the vertical one. Crucially, this perturbation still produces the exchange of Oseledets subspaces; in this case, however, we begin the analysis from the vertical direction. Let us carry out the construction assuming the strict inequality
\begin{equation}\label{est}
    2^{\left(2+\frac{\log\eta}{\log \sigma}\right)\alpha}<\eta^2.
\end{equation}
This is because the equality case in Theorem~A follows directly from the same theorem when only assuming the above strict inequality. Indeed, suppose $2^{\left(2+\frac{\log \eta}{\log \sigma}\right)\alpha} = \eta^2$. Let $\mathcal{U}$ be an open neighborhood of $A$ in $C^{\alpha}(M, \SL(2))$ and choose $r > 0$ sufficiently small so that $\eta + r < \sigma$ and the cocycle 
$$A'|_{[0;j]}=\begin{pmatrix}
(\eta+r)^{-1} & 0 \\				0 & \eta+r 
\end{pmatrix}\in \mathcal{U}.$$ 
Since $2^{\left(2+\frac{\log \eta}{\log \sigma}\right)\alpha} < (\eta + r)^2$, we may apply Theorem~A to the cocycle $A'$, which yields a cocycle $B' \in \mathcal{U}$ with zero Lyapunov exponents. This provides the desired cocycle $B'$ arbitrarily close to $A$ in the $\alpha$-H\"older topology, even when $2^{\left(2+\frac{\log \eta}{\log \sigma}\right)\alpha} = \eta^2$.   

Since Lyapunov exponents of $A$ are nonzero, for $\mu$-almost every $x \in M$, Oseledets Theorem \cite{oseledets} provides a decomposition of $\mathbb{R}^2$ into $f$-invariant subspaces called Oseledets subspaces, $\mathbb{R}^2 = E_x^s \oplus E_x^u$. Moreover, because the cocycle $A$ is defined by diagonal matrices, its Oseledets subspaces coincide $\mu$-almost everywhere with the standard vertical and horizontal bundles:
\begin{equation}\label{subos}
    \mathbb{R}^2 = H_x \oplus V_x = \mathbb{R}(1,0) \oplus \mathbb{R}(0,1).
\end{equation}
The central strategy in Bocker-Viana's approach \cite{bockerviana} is to construct a cocycle $B$, close to $A$ in the $\alpha$-Hölder topology, that exchanges the Oseledets subspaces of $A$ described in Equation \eqref{subos}. When certain conditions are met, this exchange property guarantees that Lyapunov exponents of the resulting cocycle $B$ vanish.

We apply the idea of exchanging the Oseledets subspaces and adjust the approximation procedure so that, under appropriate conditions, the cocycles $A$ become discontinuity points for the Lyapunov exponents. 
In our construction, we crucially use the difference between the expansion and contraction rates ($\sigma> \eta$) of the matrices defining the cocycle $A$.

Let us fix several key parameters. Choose $\beta \in \left[\frac{\log\eta}{\log \sigma},1\right]$ to be a sufficiently small rational number such that $2^{(2+\beta)\alpha}<\eta^2$, and let $\gamma \in [1,2)$. For each $k \in \mathbb{N}$ such that $\beta k$ is integer, we define the cylinder
\[ Z_k = [0;\underbrace{j\ldots j}_{k\text{ times}} \underbrace{r \ldots r}_{\beta k\text{ times}}] \subset M, \]
where the symbol $j$ appears $k$ times and the symbol $r$ appears $\beta k$ times. Thus, $Z_k$ has $n=k(\beta+1)$ elements. By construction, the collection $\{f^i(Z_k)\}_{i=0}^{n-1}$ is pairwise disjoint, meaning $f^i(Z_k) \cap f^j(Z_k) = \emptyset$ for $0 \le i < j \le n-1$, which in turn writing $c = (1 + \eta^{2k(\gamma-2)})^{-1/2}$, enables us to define a modified cocycle $R_k: M \to \SL(2)$ by
\begin{eqnarray}\label{b}
R_{k}(x)= \left\{
\begin{array}{ll}
\begin{pmatrix}
1 & 0 \\
\eta^{-\gamma k} & 1
\end{pmatrix}  &\mbox{if $x \in Z_k$,} \\
\\
c\begin{pmatrix}
1 & -\eta^{k(\gamma-2)} \\
\eta^{k(\gamma-2)} & 1 
\end{pmatrix}  &\mbox{if $x \in f^k(Z_k)$,}\\
\\
\begin{pmatrix}
1 & 0 \\
\eta^{k(2-\gamma)}\sigma^{-2\beta k} & 1 
\end{pmatrix}  &\mbox{if $x \in f^{n-1}(Z_k)$,}\\
\\
\text{Id}      &\mbox{otherwise.}
\end{array}\right.
\end{eqnarray}
Finally, the perturbed cocycle of interest is defined by 
\begin{equation}\label{cand}
    B_k(x) = A(x) R_k(x).
\end{equation}

Observe that cocycle $R_k$ is constant on each cylinder, and consequently both $R_k$ and $B_k$ are $\alpha$-Hölder continuous. In particular, $B_k \in C^\alpha(M, \SL(2))$ for every $k\in \N$.

We core of our proof is the selection of a specific, smaller parameter: the angle between the horizontal axis and the vector $R_k(x)e_1$ for any $x\in Z_k$. This choice permits the initial subspace exchange (from $H_x$ to $V_x$) with a reduced Holder distance between the perturbation and cocycle $A$. Thereby, yielding a stronger hypothesis for the subsequent step. 

The definition of the cocycle $A$ provides another essential ingredient. The hyperbolicity difference ($\sigma>\eta$) between the matrices forming the cocycle reduces the number of $r$-symbols needed in cylinder $Z_k$. Thus, under the stronger expansion induced by $\sigma$, the restricted cocycle $B_k|_{e_2}$ approaches the horizontal axis in fewer iterations, while maintaining a controlled Holder distance. Moreover, increasing the cylinder length enlarges the separation between points inside and outside $Z_k$, which improves the require condition.

An explicit calculation using Equations~\eqref{b} and \eqref{cand} shows that for every $x \in Z_k$,
\begin{equation}
B_{k}^{n}(x) = c\begin{pmatrix}
0 & -\eta^{(\gamma-1)k}\sigma^{\beta k} \\
\sigma^{-\beta k}(\eta^{(\gamma-3)k} + \eta^{(1-\gamma)k}) & 0 
\end{pmatrix}.
\end{equation}
From this formula, we immediately see that cocycle $B_k$ has the \emph{exchange property} on $Z_k$: for every $x \in Z_k$,
\begin{equation}\label{interchanges}
B_{k}^{n}(x)H_x = V_{f^n(x)} \quad \text{and} \quad B_{k}^{n}(x)V_x = H_{f^n(x)}.
\end{equation}
This key property will be essential for computing the Lyapunov exponents of $B_k$ in Section~3.

It remains only to prove that cocycle $B_k$ can be made arbitrarily close to $A$ in the $\alpha$-Hölder topology which is established by the following lemma.
\begin{lema}\label{casa}
For sufficiently large $k \in \mathbb{N}$ with $\beta k \in \mathbb{N}$, the $\alpha$-Hölder norm $\|A - B_k\|_{\alpha}$ becomes arbitrarily small.
\end{lema}
\begin{proof}
For the first term in Formula \eqref{norma}, by examining all cases in Formulas~\eqref{new},\eqref{b} and \eqref{cand}, we directly obtain
\[ \|A-B_k\|_{\infty} \leq \sigma \max\{\eta^{-\gamma k}, \eta^{k(2-\gamma)}, \eta^{k(2-\gamma)}\sigma^{-2\beta k}\}. \]
Since $\eta^{2-\gamma} < \eta^2$ and $\frac{\log \eta}{\log \sigma} \leq \beta$, we have $\frac{\eta^{2-\gamma}}{\sigma^{2\beta}} < 1$. Therefore, from our parameter choices of $\gamma,\beta,\eta,\sigma$, it follows that $\|A-B_k\|_{\infty} \to 0$ as $k \to \infty$.

Now we analyze the second term of Formula~\eqref{norma}. First, observe that cocycle $B_k$ modifies $A$ only on three cylinders: $Z_k$, $f^k(Z_k)$, and $f^{n-1}(Z_k)$. We consider several cases based on these cylinders. First, suppose $x$ and $y$ do not belong to any of these cylinders. From Definitions~\ref{new}, \ref{b}, and \ref{cand}, it follows that $R_k(x) = R_k(y) = \text{Id}$, and hence $B_k(x) = A(x)$ and $B_k(y) = A(y)$. Therefore, the second term vanishes, and the H\"older norm satisfies $\|A-B_k\|_{\alpha} = \|A-B_k\|_{\infty} \to 0$ as $k \to \infty$ by the argument in the previous paragraph. Next, suppose $x$ and $y$ belong to the same cylinder ($Z_k$, $f^k(Z_k)$, or $f^{n-1}(Z_k)$). By definition, $A(x) - B_k(x) = A(y) - B_k(y)$, so again the second term vanishes, and the conclusion follows as in the first case.

It only remains to analyze the cases where $x$ belongs to some cylinder ($Z_k$, $f^k(Z_k)$, or $f^{n-1}(Z_k)$) and $y$ does not. In these cases, if $A(x) \neq A(y)$, then by Definition~\ref{new}, we have $x_0 \neq y_0$, and thus $d(x, y)=1$. Applying the triangle inequality, we obtain:
\[ \frac{\|A(x) - B_k(x) - A(y) + B_k(y)\|}{d(x, y)^\alpha} \leq 2\|A - B_k\|_\infty.   \]
Therefore, the H\"older norm satisfies $\|A - B_k\|_\alpha \leq 3\|A - B_k\|_\infty \to 0$ as $k\to \infty$. For the remaining cases, we assume $A(x) = A(y)$. Using this equation and~\eqref{cand}, we can estimate the second term of~\eqref{norma} as follows:
\begin{equation}\label{ext}
    \frac{\|A(x) - B_k(x) - A(y) + B_k(y)\|}{d(x,y)^{\alpha}} \leq \frac{\|A\|_{\infty} \|R_k(x) - R_k(y)\|}{d(x,y)^{\alpha}}.
\end{equation}
For this purpose, we will consider the following cases, accounting for possibilities of $R_k$:
\begin{itemize}
    \item Case 1: $x \in Z_k$ and $y \notin Z_k$. Here we have $N(x,y) \leq k(\beta+1)$ and $d(x,y)^{-\alpha} \leq 2^{k(1+\beta)\alpha}$, hence
    \[  \frac{\|A\|_{\infty} \|R_k(x) - R_k(y)\|}{d(x,y)^{\alpha}} \leq \eta \left(\frac{2^{(1+\beta)\alpha}}{\eta^{\gamma}}\right)^k.    \]
    \item Case 2: $x \in f^{n-1}(Z_k)$ and $y \notin f^{n-1}(Z_k)$. Similarly, $N(x,y) \leq k(\beta+1)$, which implies
    \[   \frac{\|A\|_{\infty} \|R_k(x) - R_k(y)\|}{d(x,y)^{\alpha}} \leq \sigma \left(\frac{2^{(1+\beta)\alpha}\eta^{2-\gamma}}{\sigma^{2\beta}}\right)^k.    \]
    \item Case 3: $x \in f^k(Z_k)$ and $y \notin f^k(Z_k)$. By definition of cylinder $f^k(Z_k)$, $N(x,y) \leq k$ and $d(x,y)^{-\alpha} \leq 2^{k\alpha}$, yielding
    \[    \frac{\|A\|_{\infty} \|R_k(x) - R_k(y)\|}{d(x,y)^{\alpha}} \leq \sigma \left(\frac{2^{\alpha}}{\eta^{2-\gamma}}\right)^k.    \]
\end{itemize}
From these estimates and Equation~\eqref{ext}, it is clear that H\"older norm $\|A-B_k\|_\alpha \to 0$ as $k\to\infty$ provided the following three inequalities hold:
\begin{align}
2^{(1+\beta)\alpha} &< \eta^{\gamma}, \label{1} \\
2^{(1+\beta)\alpha} &< \frac{\sigma^{2\beta}}{\eta^{2-\gamma}}, \label{2} \\
2^{\alpha} &< \eta^{2-\gamma}. \label{3}
\end{align}
First, we note that inequality~\eqref{1} implies~\eqref{2}. Indeed, since $\beta \geq \frac{\log \eta}{\log \sigma}$, it follows that $\eta^{\gamma} \leq \frac{\sigma^{2\beta}}{\eta^{2-\gamma}}$, which proves inequality~\eqref{2}. Consequently, we only need to satisfy inequalities~\eqref{1} and~\eqref{3}. To optimize the parameter choice, observe that $(2-\gamma)$ is a decreasing function of $\gamma$ while $\frac{\gamma}{1+\beta}$ is an increasing one. Therefore, the optimal parameter choice occurs when $\gamma = \frac{2(1+\beta)}{2+\beta}$. With this choice, the condition
\[ 2^{(2+\beta)\alpha} < \eta^2 \]
implies both inequalities~\eqref{1} and~\eqref{3}, which completes the proof.
\end{proof}

\section{Lyapunov exponents of the perturbation}\label{casu}
This section concludes the proof of Theorem~A by showing that the perturbed cocycle $B_k$ has zero Lyapunov exponents for sufficiently large $k$. While the vanishing of Lyapunov exponents could alternatively be obtained from Theorem~8.1 of \cite{bockerviana} (since cocycle $B_k$ exchanges the Oseledets subspaces of cocycle $A$), the proof we present here is primarily based on the structure of the perturbed cocycle $B_k$ itself.

Since the construction used in the proof depends only on the constants $\eta$ and $\sigma$, and the remaining parameters have no influence in the argument, we may, for simplicity, restrict our attention to the case $l=1$. The construction for this base case $l=1$ extends naturally to all $l>1$.

In other words, we consider the following cocycle
	\begin{equation}\label{bveta}
    \begin{array}{rcl}
		\ase: \{0,1\}^{\Z} &\longrightarrow&\SL(2)\\
        \\
		(x_i)_{i \in \Z} & \longmapsto& \left\{
		\begin{array}{ll}
			\begin{pmatrix}
				\eta^{-1} & 0 \\
				0 & \eta 
			\end{pmatrix} \,\ \mbox{if $x_0=0$} \\ 
			\\
			\begin{pmatrix}
				\sigma & 0 \\
				0 & \sigma^{-1} 
			\end{pmatrix} \,\ \mbox{if $x_0=1$}.
		\end{array}\right.
        \end{array}
		\end{equation}

We henceforth denote the Bernoulli measure by $\mu_{p}=\nu_{p}^{\Z}$, where $\nu_{p}=p\delta_1+(1-p)\delta_0$ with $p \in (0,1)$.

Fix a sufficiently large $k \in \mathbb{N}$ such that $\beta k \in \mathbb{N}$. In what follows, we treat with the cylinder $Z_k=[0; 0\cdots 01\cdots 1]$ consisting of $k$ zeros followed by $\beta k$ ones. The central idea for computing the Lyapunov exponents of the system $(\{0,1\}^{\Z},f,\mu_p,B_k)$ is to restrict the problem to a subsystem based on the cylinder $Z_k$. To make this precise, let us introduce some definitions and notations. Let $Z := Z_k$ and $\mu := \mu_p$, and define $\tau_Z$ to be the first return time to $Z$ under $f$. That is, for every $x \in Z$,
\[ \tau_Z(x) = \inf\{m \geq 1 : f^m(x) \in Z\}.   \]
Furthermore, for every $m \geq 1$, let $\tau_Z^{(m)}$ denote the $m$-th return time to $Z$ under $f$. Specifically, for every $x \in Z$,
\begin{equation}\label{mreturntimes}
\tau_Z^{(m)}(x) = \sum_{i=0}^{m-1} \tau_Z\big(f^{\tau_Z^{(i)}}(x)\big),
\end{equation}
where $\tau_Z^{(m)}(x)$ counts the total number of iterations under $f$ required for $x$ to return to $Z$ exactly $m$ times. Here, we use the convention that $\tau_Z^{(0)}(x) = 0$. Let $\pi \colon \{0,1\}^{\Z} \to \{0,1\}$ be the projection onto the $0$-th coordinate. For every $m \geq 1$ and every $x \in \{0,1\}^{\Z}$, the number of occurrences of the symbol $1$ in the first $m$ coordinates of $x$ is given by
\[  S_m(x) = \sum_{i=0}^{m-1} \pi \circ f^{i}(x).  \]
Consequently, for every $x \in Z$, we can define
\[  S_{\tau_Z}(x) = \sum_{i=0}^{\tau_Z(x)-1} \pi \circ f^{i}(x).  \]
Let $\mu_Z = \frac{\mu|_Z}{\mu(Z)}$ be the normalized restriction of $\mu$ to $Z$, which is invariant under the first return map $f^{\tau_Z}$. The measure subsystem $(Z, f^{\tau_Z}, \mu_Z)$ is well known to be ergodic. Moreover, recalling Definition \ref{lll}, we obtain an induced cocycle $B^{\tau_Z}_k$ defined for $x \in Z$ by the product
\[ B^{\tau_Z}_k(x) = B_k^{\tau_Z(x)}(x) := B_k(f^{\tau_Z(x)-1}(x)) \cdots B_k(f(x))B_k(x).  \]
This gives the desired induced subsystem $(Z, f^{\tau_Z}, \mu_Z, B^{\tau_Z}_k)$. Its relevance lies in the following relation between the Lyapunov exponents:
\begin{equation}\label{relation}
\lambda_{\pm}(B^{\tau_Z}_k, \mu_Z) = \frac{\lambda_{\pm}(B_k, \mu)}{\mu(Z)}.
\end{equation}
Therefore, it suffices to show that the induced subsystem has zero Lyapunov exponents; i.e., $\lambda_{\pm}(B^{\tau_Z}_k, \mu_Z) = 0$.

Let us compute some iterates of the cocycle $B^{\tau_Z}_k$. Using Equations~\eqref{b} and \eqref{cand} along with the first return time $\tau_Z$, for every $x \in Z$ we obtain
\begin{equation}\label{matrixtau}
B_{k}^{\tau_Z}(x) = c\begin{pmatrix}
0 & -\sigma^{S_{\tau_Z}(x)}\eta^{-\tau_Z(x)+S_{\tau_Z}(x)+k\gamma} \\
\sigma^{-S_{\tau_Z}(x)}\eta^{\tau_Z(x)-S_{\tau_Z}(x)}(\eta^{-4k+k\gamma}+\eta^{-k\gamma}) & 0 
\end{pmatrix},
\end{equation}
where $c = (1+\eta^{2k(\gamma-2)})^{-1/2}$. Note that the normalization constant satisfies the identity $c^2\eta^{k\gamma}(\eta^{-4k+k\gamma}+\eta^{-k\gamma}) = 1$.

For the second return time $\tau^{(2)}_{Z}$, we have
\begin{equation}\label{matriz2}
B_{k}^{\tau^{(2)}_{Z}}(x) = 
\begin{pmatrix}
-\sigma^{R(f^{\tau_Z}(x))-R(x)}\eta^{-S(f^{\tau_Z}(x))+S(x)} & 0 \\
0 & -\sigma^{-R(f^{\tau_Z}(x))+R(x)}\eta^{S(f^{\tau_Z}(x))-S(x)}
\end{pmatrix},
\end{equation}
where we define
\begin{equation}\label{obs}
R(x) = S_{\tau_Z}(x) \quad \text{ and } \quad S(x) = \tau_Z(x) - S_{\tau_Z}(x).    
\end{equation}
By induction on $j \geq 1$, we can compute the cocycle iterates at even return times $m_j = \tau^{(2j)}_Z(x)$:
\begin{equation}\label{matriz2j}
B_{k}^{\tau^{(2j)}_{Z}}(x) = 
\begin{pmatrix}
(-1)^j\sigma^{c_j(x)}\eta^{-b_j(x)} & 0 \\
0 & (-1)^j\sigma^{-c_j(x)}\eta^{b_j(x)}
\end{pmatrix},
\end{equation}
with the exponents given by
\begin{equation}\label{cofe}
    c_j(x) = \sum_{i=1}^{j} \left(R(f^{\tau^{(2i-1)}_Z}(x)) - R(f^{\tau^{(2i-2)}_Z}(x))\right) \quad \text{ and }\quad b_j(x) = \sum_{i=1}^{j} \left(S(f^{\tau^{(2i-1)}_Z}(x)) - S(f^{\tau^{(2i-2)}_Z}(x))\right).
\end{equation}
From this representation, we obtain the norm estimate $\|B^{m_j}_k(x)\| \leq \sigma^{|c_j(x)|}\eta^{|b_j(x)|}$, which yields the following bound on the growth rate:
\begin{equation}\label{boundedlim}
\lim_{j \to \infty} \frac{1}{m_j}\log \|B^{m_j}_k(x)\| \leq \left(\lim_{j \to \infty} \frac{|c_j(x)|}{m_j}\right)\log \sigma + \left(\lim_{j \to \infty} \frac{|b_j(x)|}{m_j}\right)\log \eta.
\end{equation}
Assume temporarily that 
\begin{equation}\label{fef}
    \lim_{j\to\infty} \frac{1}{m_j}\log \|B^{m_j}_k(x)\| \leq 0.
\end{equation}
Since $B^{\tau_Z}_k$ is an $\mathrm{SL}(2,\mathbb{R})$-valued cocycle, the Furstenberg-Kesten Theorem \cite{furstenberg} implies that for $\mu_Z$-almost every $x \in Z$, the following limit exists and satisfies
\[  \lambda_+(B^{\tau_Z}_k, \mu_Z) = \lim_{m \to \infty} \frac{1}{m}\log \|B^{\tau^{(m)}_Z}_k(x)\| \geq 0. \]
Because $\{m_j\}_{j \in \mathbb{N}} = \{\tau^{(2j)}_Z\}_{j \in \mathbb{N}}$ is a subsequence of return times $\{\tau^{(m)}_Z\}_{m \in \mathbb{N}}$, the limit along this subsequence must satisfy $\lim_{j \to \infty} \frac{1}{m_j}\log \|B^{m_j}_k(x)\| = 0$. Consequently, the full limit must agree: $\lambda_+(B^{\tau_Z}_k, \mu_Z) = 0$. This establishes the vanishing of the Lyapunov exponents. 

This result, philosophically, was already expect. The definition of the induced cocycle at even return times in \eqref{matriz2j} shows that asymptotic behavior of vectors is determined by alternating expansion and contraction along invariant directions governed by coefficients $c_j(x)$ and $b_j(x)$ in \eqref{cofe}. These coefficients are sums of differences evaluated along the orbit at specific return times. Their oscillatory nature leads to cancellations, which for $\mu$-almost every point cause the average exponential growth rate--the Lyapunov exponent--to vanish. In other words, despite local expansion and contraction occur, they are balanced by the global dynamics of the induced cocycle, resulting in zero Lyapunov exponents.

Therefore, it remains only to verify the key estimate \eqref{fef}. In order to establish it, we first require the following lemma.
\begin{lema}\label{into}
    The functions $R:Z\to \R$ and $S:Z\to \R$ defined in Equation \eqref{obs} belong to $L^1(\mu_Z)$. 
\end{lema}
\begin{proof}
We begin by analyzing the function $R$. For each $m > n$, define the function $T_m \colon Z \to \mathbb{R}$ as
\[   T_m(x) := \sum_{i=n}^{m-1} \pi \circ f^{i}(x).  \]
Using the first return time $\tau_Z$, we observe that
\begin{equation}\label{ral}
R(x) = T_{\tau_Z}(x) + k\beta + 1.
\end{equation}
We note that $T_{\tau_Z}$ is a sum of random variables $X_i = \pi \circ f^{i}$ for $i > n$, which are independent and identically distributed (i.i.d.) with respect to $\mu_Z$. Furthermore, the sum $T_{\tau_Z}$ has a stopping time $\tau_Z - n$ with respect to the filtration of Borel $\sigma$-algebras on $Z$ induced by the sequence $\{X_i\}_{i > n}$. The mean and variance of $X_i$ are given by
\[  \mathbb{E}[X_i] = \int_{Z} \pi \circ f^{i} \, d\mu_Z = \int_{[i;0]} \pi \circ f^{i} \, d\mu_Z + \int_{[i;1]} \pi \circ f^{i} \, d\mu_Z = p, \]
\[  \operatorname{Var}(X_i) = \int_{Z} \left(\pi \circ f^{i} - p\right)^2 d\mu_Z = \int_{[i;0]} (0 - p)^2 d\mu_Z + \int_{[i;1]} (1 - p)^2 d\mu_Z = p(1 - p).  \]
These properties allow us to apply Wald's identity:
\[   \mathbb{E}\left[\sum_{i=1}^{\tau} X_i\right] = \mathbb{E}[X_1] \cdot \mathbb{E}[\tau],   \]
where $\{X_i\}$ are i.i.d. random variables and $\tau$ is a stopping time adapted to $\{X_i\}$. Applying this to $T_{\tau_Z}$, we obtain
\[  \int_{Z} T_{\tau_Z} \, d\mu_Z = \left(\int_{Z} \pi \circ f^n \, d\mu_Z\right) \left(\int_{Z} (\tau_Z - n) \, d\mu_Z\right) = p \left(\frac{1}{\mu(Z)} - n\right).  \]
From \eqref{ral}, we conclude that $R \in L^1(\mu_Z)$, with
\[  \int_{Z} R \, d\mu_Z = p \left(\frac{1}{\mu(Z)} - n\right) + (k\beta + 1).  \]
Finally, the integrability of $S$ follows immediately, since $S = \tau_Z - R$ is a sum of $L^1$-functions with respect to $\mu_Z$.
\end{proof}
\begin{lema}
Let $m_j = \tau_Z^{(2j)}(x)$ be the $2j$-th return time of $x$ to $Z$. Then, for $\mu_Z$-almost every $x\in Z$
\[  \limsup_{j \to \infty} \frac{1}{m_j}\log \|B^{m_j}_k(x)\| \leq 0. \]
\end{lema}
\begin{proof}
We will show that both terms on the right side of Equation \eqref{boundedlim} vanish, from which the result will follow. Let $g = f^{\tau_Z^{(2)}}$ denote the second return time map. From Equation \eqref{cofe}, we obtain
\begin{equation}\label{abirkhoff}
    \lim_{j \to \infty} \frac{1}{m_j}|c_j(x)| = \lim_{j \to \infty} \frac{j}{m_j} \cdot \lim_{j \to \infty} \frac{1}{j} \left| \sum_{i=0}^{j-1} R(g^i(f^{\tau_Z}(x))) - \sum_{i=0}^{j-1} R(g^i(x)) \right|.
\end{equation}
By Birkhoff's Ergodic Theorem and Kac's Theorem, we have
\[  \lim_{j \to \infty} \frac{m_j}{j} = \lim_{j \to \infty} \frac{1}{j} \sum_{i=0}^{2j} \tau_Z(f^{\tau_Z^{(i)}}(x)) = 2\int_Z \tau_Z \, d\mu_Z = \frac{2}{\mu(Z)}.   \]
Since $f^{\tau_Z}$ is ergodic and Bernoulli, the subsystem $(g, \mu_Z)$ is also ergodic. By Lemma \ref{into}, $R \in L^1(\mu_Z)$, and applying Birkhoff's Ergodic Theorem to \eqref{abirkhoff} yields
\[  \lim_{j \to \infty} \frac{1}{j} \sum_{i=0}^{j-1} R(g^i(f^{\tau_Z}(x))) = \int_Z R \, d\mu_Z = \lim_{j \to \infty} \frac{1}{j} \sum_{i=0}^{j-1} R(g^i(x)).  \]
Combining these results, we conclude that
\[  \lim_{j \to \infty} \frac{1}{m_j}|c_j(x)| = \frac{\mu(Z)}{2} \cdot \left| \int_Z R \, d\mu_Z - \int_Z R \, d\mu_Z \right| = 0. \]
Similarly, using Equation \eqref{cofe} to express the second term of \eqref{boundedlim} in terms of $g$, we obtain
\begin{equation}\label{bbirkhoff}
    \lim_{j \to \infty} \frac{1}{m_j}|b_j(x)| = \lim_{j \to \infty} \frac{j}{m_j} \cdot \lim_{j \to \infty}\frac{1}{j} \left| \sum_{i=0}^{j-1} S(g^i(f^{\tau_Z}(x))) - \sum_{i=0}^{j-1} S(g^i(x)) \right|.
\end{equation}
Since Lemma \ref{into} also guarantees $S \in L^1(\mu_Z)$, the same argument applied to $S$ shows that $\lim_{j} \frac{1}{m_j}|b_j(x)| = 0$.
\end{proof}

\section{Generalized Example}
As discussed in the previous section, discontinuities can arise even in simple settings, such as locally constant cocycles. This naturally leads us to ask what phenomena may occur in more general and less rigid systems.

In this section, we continue to work within the following framework: a bi-infinite sequence space over $l+1$ symbols, $l\geq 1$, and a Hölder continuous cocycle  $A:M \to \SL(2)$ of the form 
\[A(x)=\begin{pmatrix}
    a(x) & 0\\
    0 & a(x)^{-1}
\end{pmatrix},\]
where $a: M \to \R$ is a Hölder continuous function satisfying $a(x)\neq 0$ for every $x\in M$. 

Though this cocycle remains diagonal, working in this more general context requires us to impose a few additional restrictions. These nevertheless still allow much greater flexibility than the locally constant case treated earlier.

Let us consider two cylinders, $\cx$ and $\cy$, where $q_1$ and $q_2$ denote, respectively, the number of symbols in each cylinder. Thus, we write $\cx=[0; x_0 \ldots x_{q_1-1}]$ and $\cy=[0; y_0 \ldots y_{q_2-1}]$, and assume that 
\begin{eqnarray}\label{disjsets}
    f^j(\cx)\cap f^i(\cy)=\emptyset,  \,\ \,\ \,\ \forall j\in \{0, \cdots, q_1-1\} \,\ \mbox{ and }   \,\ \forall i\in \{0, \cdots, q_2-1\}.
\end{eqnarray}

\begin{exemplo}\label{cylindersm1}
    The cylinders $C_1=[0;0]$ and $C_2=[0;1]$ satisfy the above condition with $q_1=q_2=1$. Another pair $C_1=[0;010]$ and $C_2=[0;111]$ also satisfies \eqref{disjsets} now with $q_1=q_2=3$. Observe that symbols inside cylinders may repeat.
\end{exemplo}

\begin{exemplo}
    We may also consider $C_1=[0;l]$ and $C_2=[0;0 \cdots l-1]$ which satisfy condition \eqref{disjsets} with $q_1=1$ and $q_2=l$. We emphasize this definition does not restrict the values of $q_1$ and $q_2$.
\end{exemplo}
Fix two constants $\sigma\geq\eta>1$ and define the cocycle $A \colon M \to \SL(2)$ by:
\begin{equation}
\label{GC}
\begin{array}{rcl}
A: M & \longrightarrow & \SL(2) \\
\\
x & \longmapsto & \left\{
\begin{array}{lll}
\begin{pmatrix}
\eta^{-1} & 0 \\
0 & \eta
\end{pmatrix}  &\mbox{if $x \in \bigcup_{i=0}^{q_1-1}f^i(\cx)$}, \\
\\
\begin{pmatrix}
\sigma & 0 \\
0 & \sigma^{-1}
\end{pmatrix}  &\mbox{if $x \in \bigcup_{i=0}^{q_2-1}f^i(\cy)$},\\
\\
\begin{pmatrix}
				a(x) & 0 \\
				0 & a(x)^{-1} 
			\end{pmatrix} \,\ &\mbox{otherwise,}
\end{array}
\right.
\end{array}
\end{equation}
        where $a:M \to \R$ is a $\alpha$-Hölder continuous function satisfying the following restriction:
\begin{eqnarray}\label{h2}
    \int_{M\setminus C} \log a(x) \, d\mu(x) \geq 0, \quad \text{where } \quad C = \bigcup_{i=0}^{\,q_1-1} f^i(\cx) \cup \bigcup_{i=0}^{\,q_2-1} f^i(\cy).
\end{eqnarray}
This condition allows $a(x)$ to be identically $1$ outside set $C$. 

Note that $A$ is well defined because the sets $f^j\left(\cx\right)$ and $f^i(\cy)$ are disjoint for all $0\leq j\leq q_1-1$ and $0\leq i\leq q_2-1$, which is ensured by the cylinder construction.

Since $A$ coincides with the matrix formed by $a(x)$ and $a(x)^{-1}$ outside a finite union of cylinder sets, it inherits the Hölder continuity of $a(x)$, and is therefore an $\alpha$-Hölder continuous cocycle. 

\begin{obs}\label{particularcase}
        This cocycle construction generalizes the locally constant cocycle of \eqref{new}. Specifically, taking the cylinders $\cx=[0;0]$ and $\cy=[0;1]$ of Example \ref{cylindersm1} over the space $\{0,1\}^{\Z}$, yields the cocycle of \eqref{bveta}.
    \end{obs}

We also assume that the weights of the measure $\mu$ are chosen such that
\begin{eqnarray}\label{h1} 
-\log \eta \mu(\cup_{i=0}^{q_1-1}f^{i}(\cx))+\log \sigma \mu(\cup_{i=0}^{q_2-1}f^{i}(\cy))>0.  \end{eqnarray}

Under these choices, the Lyapunov exponents of $A$ are positive, which we now demonstrate. 

Suppose, by contradiction, that $\lambda_{+}(A, \mu)=0$. Since $A$ is a diagonal cocycle, we can combine Oseledets' Theorem with the Birkhoff Ergodic Theorem to compute the Lyapunov exponent as 
$$0=\lambda_{+}(A, \mu)=\int \log a(x)d\mu(x).$$ 

On the other hand, hypothesis \ref{h1} and \ref{h2} imply that the integral is positive:
\begin{eqnarray*}
    \int_{M}\log a d\mu&=& \int_{\bigcup_{i=0}^{q_1-1}f^i(\cx)}\log ad\mu+\int_{\bigcup_{i=0}^{q_2-1}f^i(\cy)}\log a d\mu+\int_{M\setminus C}\log a d\mu\\
    &=& \int_{\bigcup_{i=0}^{q_1-1}f^i(\cx)}\log \eta^{-1}d\mu+\int_{\bigcup_{i=0}^{q_2-1}f^i(\cy)}\log \sigma d\mu+\int_{M\setminus C}\log a d\mu\\
    &=&\log\sigma\mu\left(\bigcup_{i=0}^{q_2-1}f^i(\cy)\right)-\log \eta \mu\left(\bigcup_{i=0}^{q_1-1}f^i(\cx)\right)+\int_{M\setminus C}\log a d\mu\\
    &>&0.
\end{eqnarray*}

This contradiction shows that $\lambda_{+}(A, \mu)>0$. 

Observe also that cocycle $A$ is not uniform hyperbolic, since we may choose a measure $\mu$ whose associated weights satisfy $\lambda_+(A,\mu)=0$. 

We know that Backes-Brown-Butler's Theorem guarantees continuity under the $\alpha$-fiber-bunching condition. In contrast, Corollary \ref{teoa} handles the case $\max_{x \in M}\{a(x), a(x)^{-1}\}^2< 2^{\alpha}$, showing that if this domination fails, discontinuities can indeed occur, even for simple, non-locally constant diagonal cocycles over a Bernoulli shift. 

Moreover, the theorem highlights that occurrence of discontinuity can be determined solely from the local behavior of the cocycle on specific cylinders. In particular, it suffices to consider matrices that act in opposite ways on the same direction--one contracting and the other expanding the horizontal axis--to produce a discontinuity mechanism. For this reason, the proof follows the same lines as that of Theorem~A. However, since the cocycle is no longer locally constant, the argument previously used to establish the vanishing of the Lyapunov exponents cannot be applied. Instead, we prove that Lyapunov exponents are zero by using the same reasoning as in Theorem~8.1 of \cite{bockerviana}.

This example still lies far from being fiber-bunching, because, in contrast, the condition appearing in Corollary~\ref{teoa} is about the constants $\eta$ and $\sigma$, which may be much smaller than the number $\max_{x \in M}\{a(x), a(x)^{-1}\}$.

\begin{corolario}\label{teoa}
		\textit{For any $\alpha>0$ and $\eta \leq \sigma$ such that 
\begin{eqnarray}\label{condition2}
        2^{3\alpha}<\sigma^2
    \end{eqnarray}
    and
\begin{eqnarray}\label{condition31}
	2^{\left(2+\frac{\log\eta}{\log\sigma}\right)\alpha}\leq\eta^2
    \end{eqnarray}
    and for any choice of positive weights $(p_0, \cdots, p_l)$ satisfying \eqref{h1}, there exist $\alpha$-Hölder continuous cocycles $B: M \to \SL(2)$ with vanishing Lyapunov exponents which are arbitrarily close to $A$ in the $\alpha$-Hölder norm. Therefore, $A$ is a discontinuity point for Lyapunov exponents in $C^{\alpha}(M, \SL(2))$.}
	\end{corolario}
    
\section{Example in higher dimension}\label{calda}
In this section, we give a simple example of discontinuity of Lyapunov exponents in the $\alpha$-Hölder topology in the higher dimensional context, i.e., we prove Corollary \ref{dim}. The basic idea is to construct a block diagonal cocycle over a Bernoulli shift of two symbols, where each block is given by the 2-dimensional cocycles given in Equation \eqref{cla}. 

We first deal with the even dimensional case. Recall the Bernoulli shift system $(\{0,1\}^{\Z}, f, \mu_p)$ where $\{0,1\}^{\mathbb{Z}}$ is the space of bi-infinite sequences of two symbols equipped with metric $d$ given in Equation \eqref{metric}, $f \colon \{0,1\}^{\Z} \to \{0,1\}^{\Z}$ is the left shift map and $\mu_p = \nu_p^{\mathbb{Z}}$ on $\{0,1\}^{\Z}$ is the Bernoulli product measure induced by the probability measure $\nu_p = p\delta_1+(1-p)\delta_0$ on $\{0,1\}$ for $p \in (1/2,1)$.

For $m\in \N$ and $\sigma>1$, we define the cocycle $A:\{0,1\}^{\Z} \to \mathrm{SL}(2m,\mathbb{R})$ as follows
\begin{equation}\label{mem}
A(x)= 
\begin{cases}
    \begin{pmatrix}
        A_{\sigma}^{-1} & & \\
        & \ddots & \\
        & & A_{\sigma}^{-1}
    \end{pmatrix}
    & \text{if }\quad   x_0=0 \\
    \begin{pmatrix}
        A_{\sigma} & & \\
        & \ddots & \\
        & & A_{\sigma}
    \end{pmatrix}
    & \text{if }\quad   x_0=1,  \quad \text{where} \quad A_{\sigma} = \begin{pmatrix}
    \sigma & 0 \\
    0 & \sigma^{-1}
\end{pmatrix}
\end{cases}    
\end{equation}
Since $A$ is locally constant it follows that $A$ is $\alpha$-Hölder continuous. 

We also see that all Lyapunov exponents of $A$ are nonzero. Indeed, using again that the cocycle is diagonal, locally constant and the measure is Bernoulli, it follows that there are $m$ Lyapunov exponents with value $(1-2p)\log\sigma$ and $m$ exponents with value $(2p-1)\log\sigma$ where $\nu_p\{1\}=p>1/2$. In fact, the special structure of cocycle $A$ and the results of Section \ref{section2} imply that for $\mu$-almost every $x\in \{0,1\}^{\Z}$, its stable and unstable Oseledets subspaces are 
\[ E_x=\bigoplus_{j=1}^mE_{2j} \qquad \text{ and }\qquad F_x=\bigoplus_{j=1}^mE_{2j-1}, \]
where $E_j=\mathrm{Span}(e_j)$ with $j=1,\ldots, 2m$ are the canonical subspaces of $\R^{2m}$.

Similar to the 2-dimensional case, we look for perturbations $B_k$ arbitrarily close to $A$ in the $\alpha$-Hölder topology such that $B_k$ have the exchange property of Oseledets subspaces of $A$. Indeed, for $\gamma\in [1,2)$, $k\in \N$ and $n=2k+1$ we consider 
\begin{equation}\label{candas}
    Z_k = [0;\underbrace{0\ldots 0}_{k\text{ times}} \underbrace{1\ldots 1}_{k+1\text{ times}}] \subset \{0,1\}^{\Z} \quad \text{ and }\quad B_k(x) = A(x) R_k(x), 
\end{equation}
where $R_k: \{0,1\}^{\Z} \to \mathrm{SL}(2m,\R)$ is a block diagonal cocycle composed by $2\times2$ matrices similar to the block matrices given for the perturbation cocycle $R_k$ of Equation \eqref{b}. Since the restrictions of cocycles $A$ and $B_k$ to their $j$-th $2\times 2$ block matrix have the same structure as the 2-dimensional case, the results of Section \ref{section2} with $\sigma=\eta$ and $j\in\{1,\ldots,m \}$ imply that 
\begin{equation}\label{hexchange} 
B_k^n(x)(E_{2j-1})=E_{2j} \quad \text{ and }\quad  B_k^n(x)(E_{2j})=E_{2j-1}.
\end{equation}
Therefore, the cocycle $B_k$ has the exchange property of Oseledets subspaces mentioned in Section \ref{section2}:
\[ B_k^n(x)(E_x)=F_{f^n(x)} \quad \text{ and }\quad  B_k^n(x)(F_x)=E_{f^n(x)}. \]

We outline how to show that $B_k$ converges to $A$ as $k\to \infty$ in the $\alpha$-Hölder topology. For the first term in Equation \eqref{norma}, by examining all cases in Formulas ~\eqref{mem} and \eqref{candas}, we directly obtain
\[ \|A-B_k\|_{\infty} \leq \sigma \max\{\sigma^{-\gamma k},\sigma^{k(\gamma-2)}\}\longrightarrow 0 \]
as $k \to \infty$ by our choice of $\gamma$. For the second term of Equation \eqref{norma}, we deal with the same cases as proof of Lemma \ref{casa}, i.e., the different situations for which $x\in f^i(Z_k)$ and $y\in f^{j}(Z_k)$ with $i,j\in \{0,\ldots,n \}$. Thus, a similar analysis of all cases shows the quantity 
\[ \frac{\|(A-B_k)(x)-(A-B_k)(y)\|}{d(x,y)^{\alpha}} \]
is bounded above by the following expressions
\[ \sigma\left( \frac{2^{2\alpha}}{\sigma^{\gamma}}\right)^k\quad  \text{ and }\quad \sigma\left( \frac{2^{\alpha}}{\sigma^{2-\gamma}}\right)^k.  \]
Therefore, $B_k$ converges to $A$ whenever the following inequalities are satisfied
\[ 2^{2\alpha}<\sigma^{\gamma} \quad \text{ and }\quad 2^{2\alpha}<\sigma^{2(2-\gamma)}.  \]
The largest set of values of $\sigma$ fulfilling these inequalities is given by the optimal choice $\gamma=4/3$ which implies $2^{3\alpha}<\sigma^2$. Thus, $B_k$ converges to $A$ in the $\alpha$-Hölder topology if $\sigma^2>2^{3\alpha}$.  

Finally, we will show that all Lyapunov exponents of $B_k$ are zero. Suppose by contradiction there exists a nonzero Lyapunov exponent for $B_k$. Hence, there is a vector $v\neq 0$ such that
\begin{equation}\label{ado}
    \lim_{n\to \infty}\frac{1}{n}\log|B_k^n(x)v|\neq 0.
\end{equation}
First assume this limit is negative. As above, let $\tilde{P}$ be the plane generated by the coordinate basis $\{e_{2j-1},e_{2j}\}$ of $\R^{2m}$ for some $j\in \{0,\ldots, m\}$ such that $\tilde{v}=\mathrm{Proj}_{\tilde{P}}v\neq 0$. Since $B_k$ lefts $\tilde{P}$ invariant by property \eqref{hexchange}, the results of Section \ref{casu} imply that 
\[ \lim_{n\to \infty}\frac{1}{n}\log|B_k^n(x)\tilde{v}|=0 \]
because restriction of $B_k$ to the $j$-th diagonal block is exactly cocycle $B_k$ of Section \ref{section2} satisfying the exchange property of Oseledets subspaces. Since $\tilde{v}$ is a projection, negative limit in Equation \eqref{ado} shows that 
\[  0=\lim_{n\to \infty}\frac{1}{n}\log|B_k^n(x)\tilde{v}|\leq \lim_{n\to \infty}\frac{1}{n}\log|B_k^n(x)v|< 0, \]
a contradiction. If the limit in Equation \eqref{ado} is positive, an analogous reasoning using $B_k^n$ as $n\to -\infty$ provides a similar contradiction. Therefore all Lyapunov exponents of $B_k$ are zero.  

For the odd dimensional case we define an almost identical cocycle to the even case. Given $m\in \N$ and $\sigma>1$, consider the cocycle $\tilde{A}:\{0,1\}^{\Z} \to \mathrm{SL}(2m+1,\mathbb{R})$ given by 
\begin{equation}\label{impar}
\tilde{A}(x)= 
    \begin{pmatrix}
 &  &  & \vdots \\
 & A_{2m\times 2m} &  & 0 \\
 &  &  & \vdots \\
\cdots & 0 & \cdots & 1
\end{pmatrix}    
\end{equation}
where $A\in \mathrm{SL}(2m,\R)$ is the cocycle in Equation \eqref{mem}. Observe that 1 in the last entry is necessary so that $\tilde{A}\in \mathrm{SL}(2m+1,\R)$. Clearly, the cocycle $\tilde{A}$ is locally constant and is $\alpha$-Hölder continuous. Noting the special structure of $\tilde{A}$, analogous arguments to above show that $\tilde{A}$ has $m$ Lyapunov exponents with value $(1-2p)\log\sigma$, a single zero Lyapunov exponent and $m$ Lyapunov exponents with value $(2p-1)\log\sigma$. Furthermore, the Oseledets subspaces are given by
\[ E_x=\bigoplus_{j=1}^mE_{2j}, \quad F_x=\bigoplus_{j=1}^mE_{2j-1} \quad \text{ and } \quad E^0_x=E_{2m+1},  \]
where $E^0$ denotes the subspace associated to vectors with zero Lyapunov exponent. 

For every $k\in \N$, the perturbation cocycle $B_k$ will be a block diagonal cocycle formed by the perturbation cocycle of the even case plus the last diagonal entry equal to 1, like in Equation \eqref{impar}. Using the same technique it follows that
\[ B_k\longrightarrow \tilde{A} \quad  \text{ in the $\alpha$-Hölder topology if }\,\ \sigma^2>2^{3\alpha}. \]
Same argument shows that $B_k$ has all its Lyapunov exponents equal to zero. Therefore, for the odd dimensional case, $\tilde{A}$ is a point of discontinuity only for $2m$ Lyapunov exponents while for the last exponent it is a point of continuity.

\bibliographystyle{plain}
\bibliography{ref}

\end{document}